\documentclass[oneside]{amsart}
\usepackage[utf8]{inputenc}
\usepackage[english]{babel}
\usepackage{url}
\usepackage[dvips]{graphics,epsfig}
\usepackage{times}
\usepackage{amsmath,amsthm,amssymb}
\usepackage{array}
\usepackage{bbm}
\usepackage{xcolor}
\usepackage{verbatim}
\usepackage{upgreek}
\usepackage[foot]{amsaddr}
\usepackage{mathtools}

\usepackage{geometry}

\def\CCC{\mathbb{C}}

\def\FFF{\mathbb{F}}

\def\OO{\mathcal{O}}

\def\QQQ{\mathbb{Q}}

\def\UU{\mathcal{U}}

\def\ZZZ{\mathbb{Z}}

\def\seq{\subseteq}

\def\Gal{\operatorname{Gal}}
\def\Fr{\operatorname{Fr}}

\newtheorem{theorem}{Theorem}[section]

\newtheorem{lemma}[theorem]{Lemma}

\newtheorem{cor}[theorem]{Corollary}
\newtheorem{prop}[theorem]{Proposition}
\newtheorem{fact}[theorem]{Fact}

\newtheorem{conj}[theorem]{Conjecture}

\theoremstyle{definition}

\theoremstyle{remark}
\newtheorem{remark}[theorem]{Remark}
\newcommand{\m}{\mathbb }

\newcommand{\N}{\operatorname{N}}

\def\Ind{\setbox0=\hbox{$x$}\kern\wd0\hbox to 0pt{\hss$\mid$\hss} \lower.9\ht0\hbox to 0pt{\hss$\smile$\hss}\kern\wd0}
\def\Notind{\setbox0=\hbox{$x$}\kern\wd0\hbox to 0pt{\mathchardef \nn=12854\hss$\nn$\kern1.4\wd0\hss}\hbox to 0pt{\hss$\mid$\hss}\lower.9\ht0 \hbox to 0pt{\hss$\smile$\hss}\kern\wd0}

\numberwithin{equation}{section}

\newcommand{\clqedsym}{\dashv_{\text{\scriptsize claim}}}
\newcommand{\clqed}{%
  \ifmmode
    \eqno{\clqedsym}%
  \else
    \hfill$\clqedsym$%
  \fi
}

\title{Progress on the stable fields conjecture}

\date{December 3, 2026}
\author[G. Conant]{G. Conant}

\address{Department of Mathematics, Statistics, and Computer Science\\
University of Illinois Chicago}
\email{gconant@uic.edu}%

\author[J. Freitag]{J. Freitag}

\email{jfreitag@uic.edu}

\author[C. Terry]{C. Terry}

\email{caterry@uic.edu}

\author[J. Ye]{J. Ye}

\address{Department of Mathematics\\
The University of Hong Kong,
Pokfulam, Hong Kong}
\email{jinheye@hku.hk}

\thanks{GC was  supported by NSF grant DMS-2452816. JF was supported by NSF grants DMS-2452197 and DMS-2554149. CT was  supported by NSF CAREER Award DMS-2509647 and a Sloan Research Fellowship. GC and CT acknowledge OpenAI for providing complimentary access to ChatGPT through the ChatGPT for Academic Researchers program.}

\begin{document}

\begin{abstract}
We prove that every infinite stable field of positive characteristic is separably closed. We also show that  the field $\CCC(t)$ has the strict order property (hence is unstable). The latter result was independently obtained by Dembner and Franklin. Using  techniques similar to those in the proof of the first result, we also show that stable Lefschetz fields are algebraically closed, and we prove Podewski's conjecture for minimal bounded fields.
\end{abstract}

\maketitle

\section{Introduction}

A field $F$ is called \emph{stable} if its (complete) theory in the language of rings is stable. Examples of stable fields include algebraically closed fields, which are in fact $\omega$-stable (shown by Morley in 1965, cf. \cite[Example I]{Morley1965}), along with imperfect separably closed fields, which are strictly stable (shown by Macintyre, Shelah, and Wood in 1975; see \cite{Shelah1975,Wood1979}).  In the other direction, Macintyre \cite{Macintyre1971} proved in 1971 that any infinite $\omega$-stable field is algebraically closed, and this was generalized to superstable fields by Cherlin and Shelah \cite{CherlinShelah1980} in 1980. This body of work motivates the \emph{Stable Fields Conjecture}, which asserts that every infinite stable field is separably closed (recorded explicitly by Duret \cite{Duret1980} in 1980). 

Progress toward this conjecture includes (years correspond to arXiv preprints where available):
\begin{enumerate}
\item Duret \cite{Duret1980} (1980): Every stable PAC field is separably closed.
\item Poizat \cite{Poizat1983} (1983): Every infinite stable bounded field is separably closed.
\item Scanlon \cite{Scanlon1999} (1999): Every infinite stable field of positive characteristic is Artin-Schreier closed.
\item Shelah \cite{ShelahStronglyDependent} (2005): Every infinite strongly$^+$ stable field is algebraically closed.
\item Krupi\'{n}ski and Pillay \cite{KrupinskiPillay2011} (2009): Every infinite stable field of weight $1$ is separably closed.
\item Halevi and Palac\'{i}n \cite{HaleviPalacin2019} (2017): Every infinite stable field of finite dp-rank is algebraically closed.
\item Johnson, Tran, Walsberg, and Ye \cite{JTWY2024} (2020): Every large stable field is separably closed.
\item Johnson and Ye \cite{JohnsonYe2026} (2026): Every infinite strongly stable field is algebraically closed.
\end{enumerate}

The Stable Fields Conjecture has  generalizations to fields with simple theories \cite[Conjecture 1.1(2)]{PillayWalsberg2023}  or with NIP theories \cite[Conjecture 1.1]{AnscombeJahnkeNIP}. Both of these directions have been extensively studied, e.g., \cite{Ax1968,Chatzidakis1999,ChatzidakisPillay1998,
ChatzidakisVanDenDriesMacintyre1992,CKSNTP2,Duret1980,
FornasieroKaplanMatthews2026,Hrushovski2002,KaplanMatthewsWalsberg2026,
KaplanScanlonWagner2011,KimPillay1997,KimPillay1998,
MartinPizarro2006,MartinPizarroPillay2004,MartinPizarroWagner2009,
PillayPoizat1995,PillayScanlonWagner1998,PillayWalsberg2023} (in the simple case) and \cite{AnscombeJahnkeNIP,
DupontHassonKuhlmann2019,HaleviHassonStronglyDependent,
HHJStronglyDependentFields,HHJVTopologies,JahnkeHenselianExpansions,
JahnkeSimonNIP,JSWDpMinimal,JohnsonCanonicalDpMinimal,JohnsonDpFiniteVI,
JohnsonDpFiniteIA,JohnsonDpFiniteIB,JohnsonDpMinimal,
KaplanScanlonWagner2011,KKLAlmostRealClosed,ShelahStronglyDependent} (in the NIP case). For the NIP case, see also Anscombe's Bourbaki notes \cite{AnscombeBourbaki} and Johnson's recent ICM survey \cite{JohnsonICM2026}.

A notorious test case for the Stable Fields Conjecture is  $\CCC(t)$, which does not fall under any of the results described in (1) through (8) above. In 2007, Scanlon \cite{Scanlon2007slides} posed the question of whether $\CCC(t)$ is stable (see also \cite[p. 4034]{JTWY2024}). Very recently, Scanlon \cite{scanlon-ct} described several structures  definable in $\CCC(t)$, which are complicated but, nevertheless, are not known to imply instability (in fact, Scanlon shows that some of these structures are stable; see also \cite{Pillay1998}).

In the first two main results of this paper, we confirm the positive characteristic case of the Stable Fields Conjecture, and we show that in fact $\CCC(t)$ cannot  be a counterexample to the characteristic $0$ case.\footnote{This was independently obtained by Dembner and Franklin \cite{DembnerFranklin2026}; see the AI disclosure for details.}

\begin{theorem}\label{thm:SFCpos}
Every infinite stable field of positive characteristic is separably closed.
\end{theorem}

\begin{theorem}\label{thm:CT}
The  theory of $\CCC(t)$ has the strict order property witnessed by an existential formula. In particular, $\CCC(t)$ is unstable.
\end{theorem}

Let us informally describe the general idea underlying both proofs. A \emph{marked definable group} is a pair $(G,g)$ consisting of a definable group and a distinguished torsion point $g \in G$. Given a sequence of definable groups $(G_a, g_a)$ with a distinguished torsion point $g_a$, consider the relation $\varphi (X,Y)$ on pairs of marked definable groups in our sequence which holds precisely when there is a (definable) isogeny $\chi$ from $X=(G_a,g_a)$ to $Y=(G_b, g_b)$ with the property $\chi(g_a)=g_b$. Of course, it isn't surprising that this relation (in the abstract sense) has the order property. The more surprising aspect of things is that the relation $\varphi$ is definable. In $\m C(t)$, we use this relation with a sequence of marked elliptic curves with $g_a$ a $2$-torsion point and transcendental, but algebraically dependent $j$-invariants (as the associated elliptic curves are isogenous, the pairs of j-invariants in the sequence satisfy modular polynomials). In the positive characteristic case, we take $G_a=\m G_m$ to be the multiplicative group of the field under consideration, but vary the order of the torsion point. We show that in the case that the field is not separably closed, we can arrange for the relation to be definable. In both cases, we have taken pains to give the formula very specifically, and the argument that the formula is unstable can be given in purely elementary terms (for example, formalizing the proof in Lean is extremely quick and easy). We give the direct proofs of each theorem first before carefully explaining the isogeny interpretation in Section \ref{isosection}. 

The proof of Theorem \ref{thm:SFCpos} (given in Section \ref{sec:SFCpos}) consists of two fundamental components which draw from cases (1) and (3) of the Stable Fields Conjecture in the list above. First, we crucially use two consequences from Scanlon's unpublished note \cite{Scanlon1999}, namely, that if $F$ is an infinite stable field of characteristic $p>0$, then $\bar{\FFF}_p\seq F$ and $p$ does not divide the degree of any finite separable extension of $F$ (these results were extended to the NIP case by Kaplan, Scanlon, and Wagner \cite{KaplanScanlonWagner2011}). The second component bears a strong resemblance to Duret's \cite{Duret1980} proof that a non-separably closed PAC field has the independence property (generalized by Hempel \cite{Hempel2016} to obtain the $n$-independence property for all $n$). In particular, as in the Kummer case of \cite{Duret1980}, we will code combinatorial complexity  using cosets of multiplicative subgroups of the form $(F^\times)^q$ for $q\neq p$. But whereas Duret obtains the independence property using the PAC assumption, we obtain the order property using nested finite fields  and basic properties of generic elements in stable fields.

The proof of Theorem \ref{thm:CT} (given in Section \ref{sec:CT}) obtains the strict order property for the theory of $\CCC(t)$ using a formula of the form $\exists v\exists w\,\theta(x,y,v,w)$, where $\theta(x,y,v,w)$ is a  polynomial equation over $\ZZZ[t]$. We have written the argument in such a way so that it requires only elementary algebra and a basic fact in valuation theory. As a consequence, the proof appears somewhat unmotivated. However, as noted above, in 
 Section \ref{sec:iso1}, we give 
a more high-level conceptual discussion of the proof using isogenies between marked elliptic curves. Readers interested in this perspective may benefit from reviewing  this section before the proof of Theorem \ref{thm:CT}. 

In Section \ref{sec:more}, we use ideas similar to the proof of Theorem \ref{thm:SFCpos} to obtain two more results in the same vein. First, we show that any stable field of characteristic $0$, which is elementarily equivalent to an ultraproduct of fields of positive characteristic, must be algebraically closed (see Theorem \ref{thm:Lef}). Second, we show that any minimal bounded field of characteristic $0$ is algebraically closed (see Corollary \ref{cor:Pod}). This establishes a special case of Podewski's conjecture that every infinite minimal field is algebraically closed (proved in the positive characteristic case by Wagner \cite{Wagner-podewski}).

\subsection*{AI Disclosure and Concurrent Work} 

The results in this paper were obtained using ChatGPT 6.2 Astra. This followed several unsuccessful investigations by the first three authors with earlier models. A  proof of the strict order property for $\CCC(t)$ was obtained by the first three authors using  ChatGPT 6.2. This same chat then produced a proof of the Stable Fields Conjecture in positive characteristic. A more elementary proof of the $\CCC(t)$ result was then obtained after prompting the model to produce a Lean formalization. As we describe in the paper, the proof of the instability of $\m C(t)$ is, in retrospect, a small jump from the work of Scanlon \cite{scanlon-ct}, who showed that amongst a family of elliptic curves $E \rightarrow \m C$, the isogeny relation is definable in the structure $\m C(t)$. 

Independently of the above, Dembner and Franklin \cite{DembnerFranklin2026} also discovered the  proof of the strict order property for $\CCC(t)$. The first version of this paper (written by the first three authors) was submitted to the arXiv at a later time than their preprint, but during the same submission window. However, there was a one day gap between the announcements of the two preprints.

During this interim period, the fourth author used ChatGPT 6.2  to independently prove the Stable Fields Conjecture in positive characteristic in response to a prompt of whether the proof of the strict order property for $\CCC(t)$ in \cite{DembnerFranklin2026} could be applied to stable fields. This was written  in a short note circulated privately to several colleagues before the first version of this paper was announced. Following that announcement, we compared proofs and explored further ideas to produce the present version of this paper. This resulted in  simplifications to the proof of Lemma \ref{lem:divisible} (used in Theorem \ref{thm:SFCpos}) and the inclusion of Corollary \ref{cor:Pod} (Podewski's conjecture for minimal bounded fields). The proof of the strict order property for $\CCC(t)$ is unchanged from the first version. 

ChatGPT was also used for literature searches. All AI-generated proofs have been reorganized and completely rewritten by the authors. All exposition is human generated.

\section{The Stable fields conjecture in positive characteristic}\label{sec:SFCpos}

We begin with a fundamental fact about generics in stable fields, which  quotes \cite[Theorem 5.10]{Poi01}. The first part is due to Cherlin and Shelah \cite{CherlinShelah1980}, while the second is due to Poizat \cite{Poizat1983}.

\begin{fact}\label{fact:generic}
If $F$ is an infinite stable field then $(F,+)$ and $(F^\times,\cdot)$ are connected with the same generic type.
\end{fact}

The following results on stable fields of positive characteristic were proved by Scanlon in an unpublished note \cite{Scanlon1999}. The generalization to the setting of NIP fields was obtained later by Kaplan, Scanlon, and Wagner (see Corollaries 4.4 and 4.5 of \cite{KaplanScanlonWagner2011}). 

\begin{fact}\label{fact:KSW}
Let $F$ be an infinite stable field of characteristic $p>0$. Then:
\begin{enumerate}
\item $F$ has no finite separable extension of degree divisible by $p$.
\item $\bar{\FFF}_p\seq F$.
\end{enumerate}
\end{fact}

We now state an algebraic criterion for establishing that a field is separably closed. This is a well-known result in the literature on fields with model-theoretic tameness, going back at least to work of Macintyre \cite{Macintyre1971} on $\omega$-stable fields. Similar statements can be found in  \cite{CherlinShelah1980,Dupont2015,Dupont2017,DupontHassonKuhlmann2019,Duret1980,JTWY2024,JohnsonYe2026,Krupinski2015,KrupinskiPillay2011,Milliet2013,PillayWalsberg2023,Poi01,Sinclair2018,Wagner1998,Yoneda2006}.   The result has direct precedents in the algebra literature, including \cite[Definition 5.17]{Moerman2024} and the proof of Theorem 12 in \cite[6.\S5]{Ramanathan1954}. More general cohomological statements can be found in \cite[\S3]{Douady1960} as well as  \cite[I.\S3.3, Corollary 1 and I.\S4.1, Proposition 21]{Serre2002} (see also \cite[Proposition 2.2]{IzquierdoLucchiniArteche2025}). For the sake of completeness, we  include the argument.

\begin{lemma}\label{lem:test}
Let $F$ be a field such that $\textnormal{char}(F)$ does not divide the degree of any finite Galois extension of $F$. Suppose that for any prime $q\neq\textnormal{char}(F)$, if $K$ is  a finite separable extension  of $F$, then $K^q=K$. Then $F$ is separably closed.
\end{lemma}
\begin{proof}
Suppose $F$ is not separably closed. Let $q$ be the smallest prime dividing the degree of some finite Galois extension of $F$. Then $q\neq \text{char}(F)$ by assumption. Let $\zeta$ be a primitive $q^{\text{th}}$ root of unity. Then $F(\zeta)$ is a Galois extension of $F$ with degree dividing $q-1$, and hence its degree is $1$ by minimality of $q$. So $\zeta\in F$. 

Now let $L$ be a finite Galois extension of $F$ with degree divisible by $q$.   Let $K$ be the fixed field of a subgroup of $\Gal(L/F)$ of order $q$. So $\Gal(L/K)$ is cyclic of order $q$. Since $\zeta\in F\seq K$, we may apply \cite[\S14.7, Proposition 37]{DummitFoote2004} (see also \cite[Tag 09DX]{stacks-project}) to conclude that $L=K(c)$ for some $c\in L$ such that $c^q\in K$. Since $K^q=K$ by assumption, we have $c^q\in K^q$ and hence $c=\zeta^t d$ for some $d\in K$ and $0\leq t<q$. But now $c\in K$ since $\zeta,d\in K$, contradicting that $L=K(c)$ is a proper extension.  
\end{proof}

The main work done here will be to prove that any infinite stable field of positive characteristic satisfies the conditions of the previous lemma. In order to do this, we will need the following simple observation. We were unable to find this exact statement in the literature, but it is close to parts of previous work of Duret and others.

\begin{prop}\label{prop:pick}
Suppose $F$ is an infinite stable field and $H$ is a proper definable subgroup of $F^\times$. Then for every finite $C\subset F$, there is some $b\in F\backslash C$ such that the map $c\mapsto (b-c)H$ is injective on $C$.
\end{prop}
\begin{proof}
Let $K\succ F$ be a sufficiently saturated elementary extension and let $b^*\in K$ realize the unique (additive and multiplicative) generic type over $F$. For any distinct $c,c'\in F$, we have 
\[
\frac{b^*-c}{b^*-c'}=1+\frac{c'-c}{b^*-c'},
\]
and hence $(b^*-c)/(b^*-c')$ is also generic over $F$ since genericity is preserved by inversion and by additive (resp., multiplicative)  translation by elements of $F$ (resp., $F^\times$). Since $F^\times$ is connected, $H$ has infinite index and thus no generic element lies in $H(K)$. Therefore $(b^*-c)/(b^*-c')\not\in H(K)$ for all distinct $c,c'\in F$. The proposition now follows by elementarity. 
\end{proof}

We can now prove the main lemma.

\begin{lemma}\label{lem:divisible}
Suppose $F$ is an infinite stable field of characteristic $p>0$. Then for any prime $q\neq p$, $F^q=F$.
\end{lemma}
\begin{proof}
Fix a prime $q\neq p$ and, toward a contradiction, suppose $F^q\neq F$. Then $H\coloneqq (F^\times)^q\neq F^\times$. Note that $H$ is a definable  subgroup of $F^\times$. Let $\varphi(x,y;z)$ be a formula (in the language of rings) defining
\[
(x-z)H=(y-z)H.
\]
We will show that $\varphi(x,y;z)$ has the order property, contradicting stability of $F$. 

By Fact \ref{fact:KSW}(2), $\bar{\FFF}_p\seq F$. Given $n\geq 1$, let $p_n=p^{(q-1)n}$ and note that $p_n\equiv_q 1$. Now fix some arbitrary $n\geq 1$ and, using Proposition \ref{prop:pick}, choose $b\in F\backslash \FFF_{p_n}$ such that $c\mapsto (b-c)H$ is injective for $c\in \FFF_{p_n}$. Consider some $i\geq 1$ dividing $n$. Given $c\in \FFF_{p_n}$, we have 
\[
(b^{p_i}-c)H=(b-c^{1/p_i})^{p_i}H=(b-c^{1/p_i})H,
\]
where the second equality uses $p_i\equiv_q 1$. Therefore, for $c\in\FFF_{p_n}$, 
\begin{multline*}
F\models\varphi(b^{p_i},b;c)\quad\Leftrightarrow\quad (b^{p_i}-c)H=(b-c)H\quad\Leftrightarrow\quad (b-c^{1/p_i})H=(b-c)H\\
\quad\Leftrightarrow\quad c^{1/p_i}=c\quad\Leftrightarrow\quad c=c^{p_i} \quad\Leftrightarrow\quad c\in\FFF_{p^i}.
\end{multline*}
In other words, the trace of $\varphi(b^{p_i},b;z)$ on $\FFF_{p_n}$ is precisely $\FFF_{p_i}$. 

We now deduce the order property for $\varphi(x,y;z)$. In particular, given $m\geq 1$, apply the above with $n=2^m$ and, for $1\leq j\leq m$ choose $c_j\in\FFF_{p_{2^j}}\backslash \FFF_{p_{2^{j-1}}}$. Then for $i,j\leq m$, 
\[
F\models \varphi(b^{p_{2^i}},b;c_j)\quad\Leftrightarrow\quad c_j\in\FFF_{p_{2^i}}\quad\Leftrightarrow\quad i\geq j.\qedhere
\] 
\end{proof}

Finally, we deduce the positive characteristic case of the Stable Fields Conjecture.

\begin{proof}[Proof of Theorem \ref{thm:SFCpos}]
Let $F$ be an infinite stable field of characteristic $p>0$. We want to prove that $F$ is separably closed.
Toward applying Lemma \ref{lem:test}, we first note that $p$ does not divide the degree of any finite Galois extension of $F$ by Fact \ref{fact:KSW}(1). Now fix a prime $q\neq p$ and a finite separable extension $K$ of $F$. Then $K$ is interpretable in $F$ (see, e.g., \cite[Proposition 3.2]{scanlon-ct}), and hence is still stable. So $K^q=K$ by Lemma \ref{lem:divisible}. Since $q$ and $K$ were arbitrary, $F$ is separably closed by Lemma \ref{lem:test}. 
\end{proof}

\section{$\CCC(t)$ is unstable}\label{sec:CT}

We will use:
\begin{itemize}
\item[$\ast$] $a,b,c,\ldots$ for elements of $\CCC$, 
\item[$\ast$] $p,q,r,\ldots$ for elements of $\CCC[t]$, and 
\item[$\ast$] $f,g,h,\ldots$ for elements of $\CCC(t)$. 
\end{itemize}

Let $\theta(x,y,v,w)$ be the following formula in the language of rings (with parameters in $\ZZZ[t]$)
\[
 ((t+x)v^2-y)^2=(t^2-4)w^2+4.
\]
Let $\varphi(x,y)$ be $\exists v\exists w\,\theta(x,y,v,w)$.

\begin{lemma}\label{lem:one}
Fix $a,b\in\CCC$ with $a\neq \pm 2$. Then $\CCC(t)\models\varphi(a,b)$ if and only if $\varphi(\CCC(t),a)\seq\varphi(\CCC(t),b)$.
\end{lemma}
\begin{proof}
Assume first that $\varphi(\CCC(t),a)\seq\varphi(\CCC(t),b)$. Note that $\CCC(t)\models\theta(a,a,1,1)$, and hence $a\in \varphi(\CCC(t),a)$. So $a\in\varphi(\CCC(t),b)$ by assumption, i.e., $\CCC(t)\models\varphi(a,b)$.

Conversely, assume $\CCC(t)\models\varphi(a,b)$ and fix some $k\in\CCC(t)$ with $\CCC(t)\models\varphi(k,a)$. We want to show $\CCC(t)\models\varphi(k,b)$.
Choose $f_a,g_a,f_b,g_b\in\CCC(t)$ such that $\CCC(t)\models\theta(a,b,f_b,g_b)$ and $\CCC(t)\models\theta(k,a,f_a,g_a)$. Define $h=(t+k)f_a^2-a$. Since $\CCC(t)\models\theta(k,a,f_a,g_a)$ we have 
\begin{equation}\label{eq:htheta1}
h^2=(t^2-4)g_a^2+4.
\end{equation}
We now consider two cases depending on whether $h$ is in $\CCC$.

Suppose first that $h\not\in\CCC$. Let $\sigma\colon\CCC(t)\to \CCC(t)$ be the field embedding which sends $t$ to $h$ and fixes $\CCC$ pointwise. Since $\CCC(t)\models\theta(a,b,f_b,g_b)$, we have
\begin{equation}\label{eq:abfg}
((t+a)f_b^2-b)^2=(t^2-4)g_b^2+4.
\end{equation}
Applying $\sigma$ to (\ref{eq:abfg}) yields
\begin{equation}\label{eq:sigabfg}
((h+a)\sigma(f_b)^2-b)^2=(h^2-4)\sigma(g_b)^2+4.
\end{equation}
Since $h+a=(t+k)f^2_a$ (by definition of $h$) and $h^2-4=(t^2-4)g_a^2$ (by (\ref{eq:htheta1})), equation (\ref{eq:sigabfg}) yields
\begin{equation}\label{eq:tkcomp}
((t+k)(f_a\sigma(f_b))^2-b)^2=(t^2-4)(g_a\sigma(g_b))^2+4.
\end{equation}
From (\ref{eq:tkcomp}), we obtain $\CCC(t)\models\theta(k,b,f_a\sigma(f_b),g_a\sigma(g_b))$. So $\CCC(t)\models\varphi(k,b)$, as desired. 

Finally, suppose that $h\in\CCC$. By (\ref{eq:htheta1}), $(h^2-4)/(t^2-4)$ is a square in $\CCC(t)$, and hence must be $0$. So $h^2=4$. Since $a\neq\pm 2$, we have $a+h\neq 0$. Choose $c\in\CCC$ such that $c^2=(b+2)/(a+h)$. Then
\begin{equation}\label{eq:tkc}
b+2=(a+h)c^2=(t+k)(f_ac)^2,
\end{equation}
where the second equality uses the definition of $h$. Now (\ref{eq:tkc}) yields $((t+k)(f_ac)^2-b)^2=4$, i.e., $\CCC(t)\models\theta(k,b,f_ac,0)$. So $\CCC(t)\models\varphi(k,b)$, as desired. 
\end{proof}

\begin{lemma}\label{lem:two}
Suppose $a\in\CCC$ is transcendental. Then there is some transcendental $b\in\CCC$ such that 
\[
\CCC(t)\models\varphi(a,b)\wedge\neg\varphi(b,a).
\]
\end{lemma}
\begin{proof}
Let $c\in\CCC$ be a solution to $x^2+ax+1=0$ and define $b=2-4c^2$. Note that $b$ and $c$ are transcendental.  Let $d\in\CCC$ be a square root of $c$. Then one can directly verify 
\[
\CCC(t)\models\theta\left(a,b,\frac{d(t-2c)}{t+a},\frac{ct-2}{t+a}\right).
\]
So $\CCC(t)\models\varphi(a,b)$. 

It remains to show that $\CCC(t)\models\neg\varphi(b,a)$. First, let $R=\QQQ[c]_{(c)}$ and note that $R$ is a local subring of $\CCC$ with maximal ideal $cR$. By \cite[Theorem 10.2]{Matsumura1986} (see also \cite[Tag~00IA]{stacks-project}), there is a valuation ring $\OO$ of $\CCC$ such that $R\seq\OO$ and $\mathfrak{m}_{\OO}\cap R=cR$. Let $\boldsymbol{k}$ be the residue field of $\OO$. Given $p\in\OO[t]$, let $\bar{p}\in\boldsymbol{k}[t]$ be obtained by reducing the coefficients of $p$ modulo $\mathfrak{m}_{\OO}$. Note that $\bar{c}=0$ since $c\in\mathfrak{m}_{\OO}$. Also, this reduction map induces an embedding of $\QQQ$  in $\boldsymbol{k}$ since $\mathfrak{m}_{\OO}\cap \QQQ=\{0\}$.

Now, toward a contradiction, suppose $\CCC(t)\models\varphi(b,a)$ and fix $f,g\in\CCC(t)$ such that $\CCC(t)\models\theta(b,a,f,g)$. Set $h=(t+b)f^2-a$. Then from $\CCC(t)\models \theta(b,a,f,g)$, we obtain
\begin{equation}\label{eq:htf4}
h^2=(t^2-4)g^2+4.
\end{equation}
Also, using $d^2=c$, $c^2+1=-ca$, $b=2-4c^2$, and the definition of $h$, one can directly verify
\begin{equation}\label{eq:chtdr}
ch=(t+2-4c^2)(df)^2+c^2+1.
\end{equation}

\noindent\textit{Claim.} There are  polynomials $p_h,p_g,q,r,s\in\OO[t]$ such that 
\[
h=\frac{p_h}{q},\quad g=\frac{p_g}{q},\quad df=\frac{r}{s},\quad (\bar{p}_h,\bar{p}_g,\bar{q})\neq (0,0,0),\quad\text{and}\quad (\bar{r},\bar{s})\neq (0,0).
\]

\noindent\textit{Proof.} This is a standard consequence of the fact that $\OO$ is a valuation ring of $\CCC$. More generally, consider finitely many $f_1,\ldots,f_n\in\CCC(t)$ and write $f_i=p_i/q$ for some $p_1,\ldots,p_n\in\CCC[t]$ and common denominator $q\in \CCC[t]$. Among all 
 nonzero coefficients appearing in $p_1,\ldots,p_n,q$, we may choose a coefficient $\alpha$ of least valuation, i.e., with $\beta/\alpha\in\OO$ for all coefficients $\beta$. Now replace $p_1,\ldots,p_n,q$ with $\alpha^{-1}p_1,\ldots,\alpha^{-1}p_n,\alpha^{-1}q$. Then $p_1,\ldots,p_n,q\in\OO[t]$ and we still have $f_i=p_i/q$ for all $i$. Moreover, $(\bar{p}_1,\ldots,\bar{p}_n,\bar{q})\neq (0,\ldots,0,0)$ since $1$ appears as a coefficient in at least one of $p_1,\ldots,p_n,q$. 
 
 The claim now follows by applying this general fact separately to $h,g$ and to $df$.\clqed\medskip

Let $p_h,p_g,q,r,s\in\OO[t]$ be as in the claim. Clearing denominators in (\ref{eq:htf4}) and (\ref{eq:chtdr}) yields
\begin{align}
p^2_h &= (t^2-4)p_g^2+4q^2 \label{eq:phtf4}\\
cp_hs^2 &=q((t+2-4c^2)r^2+(c^2+1)s^2).\label{eq:pchtdr}
\end{align}
Reducing modulo $\mathfrak m_{\OO}$ then gives the following identities in $\boldsymbol{k}[t]$ (recall $\bar{c}=0$): 
\begin{align}
\bar{p}^2_h &= (t^2-4)\bar{p}_g^2+4\bar{q}^2 \label{eq:rphtf4}\\
0 &=\bar{q}((t+2)\bar{r}^2+\bar{s}^2).\label{eq:rpchtdr}
\end{align}
We now consider two cases depending on whether $\bar{q}=0$. 

First assume $\bar{q}=0$. Then $\bar{p}_g\neq 0$ since otherwise $(\bar{p}_h,\bar{p}_g,\bar{q})=(0,0,0)$ by (\ref{eq:rphtf4}). But now (\ref{eq:rphtf4}) implies that $t^2-4$ is a square in $\boldsymbol{k}(t)$, which is a contradiction. 

Finally, suppose $\bar{q}\neq 0$. Then $(t+2)\bar{r}^2+\bar{s}^2=0$ by (\ref{eq:rpchtdr}). But this implies $(\bar{r},\bar{s})=(0,0)$, since otherwise $-(t+2)$ is a square in $\boldsymbol{k}(t)$. 
\end{proof}

We now use the two previous lemmas to deduce the following sharpening of Theorem \ref{thm:CT}.

\begin{theorem}
There is a sequence $(a_i)_{i<\omega}$ of transcendental elements of $\CCC$ such that for all $i<\omega$, $\varphi(\CCC(t),a_i)\subsetneq\varphi(\CCC(t),a_{i+1})$.  Consequently, the complete theory of $\CCC(t)$ has the strict order property witnessed by the existential formula $\varphi(x,y)$.
\end{theorem}
\begin{proof}
By iterating Lemma \ref{lem:two}, we obtain a sequence $(a_i)_{i<\omega}$ of transcendental elements of $\CCC$ such that for all $i<\omega$, $\CCC(t)\models\varphi(a_i,a_{i+1})\wedge\neg\varphi(a_{i+1},a_i)$. Now Lemma \ref{lem:one} implies that $\varphi(\CCC(t),a_i)\subsetneq \varphi(\CCC(t),a_{i+1})$ for all $i<\omega$. 
\end{proof}
    
\begin{remark}
The proof implies that for $i,j<\omega$, 
\[
\CCC(t)\models\varphi(a_i,a_j)\quad\text{if and only if}\quad i\leq j.
\]
Thus $\varphi(x,y)$ induces a linear order on the sequence $(a_i)_{i<\omega}$ itself.
\end{remark}

\section{Isogenies} \label{isosection}
\subsection{Marked Isogenies in $\m C(t)$}\label{sec:iso1}

In this section, a \emph{marked elliptic curve} is a pair $(E,T)$ consisting of an elliptic curve $E$ over $\CCC$ and a chosen point $T$ of order two. An isogeny $\psi\colon E\to E'$ from $(E,T)$ to $(E',T')$ is \emph{marked} if $\psi(T)=T'$. This is the $n=2$ case of the classical moduli problem of elliptic curves equipped with a point of order $n$ \cite[Chapter~IV, Construction~4.8, p.~224]{deligne-rapoport}.

For $a\in \m C$ with $a\neq\pm 2$, let $E_a$ be the elliptic curve with affine equation $$y^2=x(x^2+ax+1).$$ Note that this cubic has distinct roots.  Write $O_a$ for the origin at infinity and $T_a=(0,0)$ for the marked point. The negation formula gives $-T_a=T_a$, so $T_a$ has order two \cite[Chapter~III,  2.3(a)]{silverman-aec}. Note that on this elliptic curve $-(x,y)= (x,-y)$. 

Consider the relation between $a$ and $b$ given by $E_a$ isogenous to $E_b$. The relation (given that the classes are countably infinite) is not definable in the field $\m C$, and in this sense seems like a good candidate for the wild behaviour of any structure where the relation is definable. This is related to the approach taken by Scanlon \cite{scanlon-ct}. Scanlon's results are positive in the sense that he shows a number of related expansions of the field $\m C$ which are $\m C(t)$-definable are stable and decidable. 

The problem from our perspective is that isogeny is symmetric, and thus won't be useful for getting an instance of the order property directly for the relation itself. The basic idea here is that replacing isogeny with isogeny of marked elliptic curves kills the symmetry and yields an instance of the order property. We give a sketch below that this is indeed a definable relation in $\m C(t)$, where the defining formula is precisely the one used to obtain the strict order property in Section \ref{sec:CT}, namely,
\[
\varphi(x,y)\coloneqq \exists v\exists w\big(((t+x)v^2-y)^2=(t^2-4)w^2+4\big).
\]

\begin{lemma}\label{lem:quartic-marked-isogeny}
Fix $a,b\in\CCC$ with $a,b\neq\pm 2$. Then
$\CCC(t)\models\varphi(a,b)$ if and only if there is an isogeny
$\psi\colon E_a\to E_b$ such that $\psi(T_a)=T_b$.
\end{lemma}

\begin{proof}
In the function field $\CCC(E_a)$, identify $t$ with the nonconstant
function $x+x^{-1}$. The equation of $E_a$ then gives
$$\left(\frac{y}{x}\right)^2=t+a,\qquad (x-x^{-1})^2=t^2-4.$$
For a general point $(x,y)$, the line through $(x,y)$ and $T_a$ has slope $y/x$. The two nonzero abscissae of intersection with $E_a$ have product $1$, so the third point has coordinates $(x^{-1},y/x^2)$. The chord-and-tangent law therefore gives \cite[Chapter~III, 2.1, 2.3(c)]{silverman-aec} $$(x,y)+T_a=(x^{-1},-y/x^2).$$ In particular, translation by $T_a$ fixes $t$ and negates both $y/x$ and $x-x^{-1}$.

Suppose first that $\m C (t)\models\varphi(a,b)$. Choose $f,g\in \m C(t)$ such that $\m C(t)\models\theta(a,b,f,g)$. Define $h=(t+a)f^2-b$, as in Section~3, so that
$$h^2=(t^2-4)g^2+4.$$
Set $$X=\frac{h+g(x-x^{-1})}{2},\qquad Y=Xf\frac{y}{x}.$$ The quartic shows that the product of $X$ and $(h-g(x-x^{-1}))/2$ is $1$. Hence $$X+X^{-1}=h,\qquad Y^2=X^2(t+a)f^2=X^2(h+b)=X(X^2+bX+1).$$
Thus $(X,Y)$ defines a rational map $F\colon E_a\dashrightarrow E_b$. Note the poles of this rational map are exactly at the points of $E_a$ mapping to $O_b.$ Since $E_a$ is smooth and $E_b$ is projective, this map extends uniquely to a morphism \cite[Chapter~II, Proposition~2.1]{silverman-aec}.

Translation by $T_a$ fixes $f,g,h$ and replaces
$g(x-x^{-1})$ and $fy/x$ by their negatives. It therefore sends $(X,Y)$ to $(X^{-1},-Y/X^2)$. The same addition formula on $E_b$ identifies this pair with the result of translation by $T_b$.  Now, because the two morphisms agree on a dense open, it follows that $$F(P+T_a)=F(P)+T_b.$$

Set $\psi(P)=F(P)-F(O_a)$. Translation on $E_b$ is a morphism, so $\psi$ is a morphism taking $O_a$ to $O_b$ and therefore a group homomorphism (a morphism of curves between two elliptic curves which fixes the identity is a group homomorphism). Evaluating the translation identity at $O_a$ gives $\psi(T_a)=T_b$. Since $T_b\neq O_b$, $\psi$ is nonconstant. It is thus a surjective homomorphism with finite kernel, as required.

Now we consider the converse. Suppose that $\psi\colon E_a\to E_b$ is an isogeny of elliptic curves with $\psi(T_a)=T_b$. On $E_c$, for $c=a,b$, let
 $$S_c=\frac{y}{x},\qquad W_c=x-x^{-1},$$
using the affine coordinates of $E_c$. Negation acts on $(S_c,W_c)$ as $(-S_c,W_c)$, and translation by $T_c$ acts as $(-S_c,-W_c)$. On $E_a$ these two automorphisms fix $t$ and generate any change of signs to the affine coordinates $(S_a,W_a)$. Also
\[
 \m C (E_a)=\m C (t)(S_a,W_a),
\]
since $x=(t+W_a)/2$ and $y=xS_a$. The identities $S_a^2=t+a$ and $W_a^2=t^2-4$ give degree at most four over $\m C(t)$. It follows that the common fixed field of the two automorphisms is $\m C(t)$.

The isogeny $\psi$ commutes with negation and intertwines translation by $T_a$ with translation by $T_b$. Therefore the functions $$f=\frac{\psi^*S_b}{S_a},\qquad g=\frac{\psi^*W_b}{W_a}$$ are fixed by both automorphisms of $E_a$, so $f,g\in\m C (t)$. Pulling back the identity $W_b^2=(S_b^2-b)^2-4$ gives
\[
 ((t+a)f^2-b)^2=(t^2-4)g^2+4.
\]
Thus $\m C (t)\models\theta(a,b,f,g)$, and hence
$\m C (t)\models\varphi(a,b)$.
\end{proof}

\subsection{Marked isogenies in the positive characteristic proof}\label{sec:iso2}
In Subsection~\ref{sec:iso1}, the formula $\varphi$ detects when there is an isogeny of elliptic curves carrying one specified marked point to another. This relation ends up encoding the order property. The techniques used in our positive characteristic proof of the stable fields conjecture can be viewed somewhat analogously, as we now explain. In this case, we work with a fixed group, namely $\m G_{\mathrm m}$, but vary the marked point. 

Let $F$ be an infinite stable field of characteristic $p>0$. By Fact~\ref{fact:KSW}(2), we have $\bar{\FFF}_p\seq F$. Fix a prime $q\neq p$. Set $Q=p^{q-1}$ and let $\sigma(x)=x^Q$ on $F$. For $n\geq 1$, set
\[
 k_n=\FFF_{Q^{p^n}},\qquad w_n=Q^{p^n}-1,
\]
and choose a generator $c_n$ of $k_n^\times$.

The self-isogenies of $\m G_{\mathrm m}$ are the power maps $[w]\colon x\mapsto x^w$ for nonzero integers $w$. To see this, recall that the invertible regular functions on $\m G_{\mathrm m}$ are the Laurent monomials $ax^w$ with $a\ne 0$ and $w\in\m Z$ (see, e.g., \cite[p.~90]{LamFirstCourse}). Moreover, for a self-isogeny, we must have $a=1$ (since the identity is preserved) and $w\neq 0$ (by surjectivity as a morphism of algebraic groups). For each $c\in k_n^\times$, by choice of $c_n$ there is a positive integer $w$ such that $[w](c_n)=c$. In particular, a marked isogeny from $(\m G_{\mathrm m},c_i)$ to $(\m G_{\mathrm m},c_j)$ exists exactly when $i\geq j$. 

The set of possible images of $c_n$ is also the kernel of the isogeny $[w_n]$. For $c\in
F^\times$ we have $[w_n](c)=\sigma^{p^n}(c) c^{-1}$, and thus 
\[
 \ker([w_n])(F)=k_n^\times.
\]
Therefore, we have  a self-isogeny of $\m G_m$ sending $c_n$ to $c$ if and only if $\sigma^{p^n}(c)=c$. This property of the power of the Frobenius map is precisely what we used in Lemma \ref{lem:divisible}. 

Suppose now that $H=(F^\times)^q$ is a proper subgroup of $F^\times$. Let $\varphi(x,y;z)$ be the formula $(x-z)H=(y-z)H$. Given $n\geq 1$, use Proposition \ref{prop:pick} to choose $b\in F\backslash k_n$  so that the cosets $(b-c)H$, for $c\in k_n$, are pairwise distinct. Then, using the argument in the proof of Lemma \ref{lem:divisible}, we conclude that $F\models \varphi(\sigma^{p^i}(b),b;c)$ if and only if $c\in k^\times_i$. For $1\leq i,j\leq n$, since $k_i^\times=\langle c_i\rangle$, we have
\[
F\models\varphi(\sigma^{p^i}(b),b;c_j) \quad\Leftrightarrow\quad c_i^w=c_j \text{ for some } w>0.
\]
In other words, $F\models\varphi(\sigma^{p^i}(b),b;c_j)$ if and only if there is a
self-isogeny of $\m G_{\mathrm m}$ over $F$ sending $c_i$ to $c_j$.
Since $k_j^\times=\langle c_j\rangle$, this happens exactly when
$k_j^\times\subseteq k_i^\times$, or equivalently, when $j\leq i$.

\section{Further examples with power class maps}\label{sec:more}

In this section, we prove two more results which showcase the same ideas behind the proof of Theorem \ref{thm:SFCpos}, namely, the interaction between power classes and standard algebraic operators. First we show that stable Lefschetz fields are algebraically closed using  power classes, Frobenius, and Duret's proof of the independence property for infinite pseudofinite fields \cite{Duret1980}. We then use power classes and the norm operator  to prove that an infinite bounded minimal field is algebraically closed (a special case of Podewski's conjecture).

\subsection{Stable Lefschetz fields}

Following Anscombe and Fehm \cite{AnscombeFehmFragments}, we define an \emph{$\mathbf{F}_{>0}$-field} to be a field $F$ for which every  sentence true in $F$ (in the language of rings) is true in some field of positive characteristic. Equivalently, an $\mathbf{F}_{>0}$-field is one   elementarily equivalent to an ultraproduct of positive characteristic fields. Such ultraproducts of characteristic $0$ are called \emph{Lefschetz fields} by Aschenbrenner and Schoutens  \cite{AschenbrennerSchoutens2007}.

\begin{theorem}\label{thm:Lef}
Every  stable $\mathbf{F}_{>0}$-field of characteristic $0$ is algebraically closed.
\end{theorem}
\begin{proof}
Since stability and algebraic closedness are properties of the theory, it suffices to consider a stable field $F$ of the form  $\prod_{\UU}F_i$  with $\text{char}(F)=0$ and $\text{char}(F_i)=p_i>0$. By \cite[Lemma 3.23$(ii)$]{Nguyen2024}, any finite extension $K$ of $F$ is also an ultraproduct of fields of positive characteristic (up to isomorphism). Such a $K$ is also stable (as in the proof of Theorem \ref{thm:SFCpos}). Therefore, to show that $F$ is algebraically closed, it suffices by Lemma \ref{lem:test} to show that $F^q=F$ for any prime $q$. We will do this with an argument similar to the proof of Lemma \ref{lem:divisible}.

Fix a prime $q$ and for a contradiction, suppose $F^q\neq F$. So $H\coloneqq (F^\times)^q$ is a proper multiplicative subgroup.
Let $\Fr_i\colon F_i\to F_i$ be the Frobenius map and let $\Fr=\prod_{\UU}\Fr_i$. Set $\sigma=\Fr^{q-1}$ and note that $\sigma(a)H=aH$ for all $a\in F^\times$ since $p^{q-1}_i\equiv_{q} 1$ for $\UU$-almost all $i$ (note that $p_i\to_{\UU}\infty$). Let $k=\textnormal{Fix}(\sigma)=\prod_{\UU}\textnormal{Fix}(\text{Fr}_i^{q-1})$. Then  $k$ is pseudofinite and infinite since $\FFF_{p_i}\seq  \text{Fix}(\text{Fr}_i^{q-1})$ and $p_i\to_{\UU}\infty$. By \cite[Lemme 6.2]{Duret1980}, the relation ``$x+y$ is a square'' has the independence property in $k$. We will lift this to witness the independence property in $F$. Precedents for this idea can be found in \cite[Example 7.7(2)]{Chernikov2014} and \cite[pp. 2998-2999]{CH99} (the latter  involves a pseudofinite fixed field of the automorphism of a difference field).

Fix $m\geq 1$. By the above, we may choose $a_1,\ldots,a_m\in k$ and $b_S\in k$ for $S\seq \{1,\ldots,m\}$, such that $a_i+b_S$ is a square in $k$ if and only if $i\in S$. Let $D$ be the set of square roots in $F$ of all sums $a_i+b_S$ for $i=1,\ldots,m$ and $S\seq \{1,\ldots,m\}$. Since $\sigma$ fixes these sums, it induces a permutation of $D$. By Proposition \ref{prop:pick}, there is some $c\in F\backslash D$ such that the cosets $(c-d)H$ for $d\in D$ are pairwise distinct.

Given $d\in D$, if $e$ is the unique element of $D$ with $\sigma(e)=d$, then we have 
\[
(\sigma(c)-d)H=\sigma(c-e)H=(c-e)H,
\]
and so
\[
(\sigma(c)-d)H=(c-d)H\quad\Leftrightarrow\quad e=d \quad\Leftrightarrow\quad d=\sigma(d)\quad\Leftrightarrow\quad d\in k.
\]
Now let $\varphi(x;y,u,v)$ be the formula
\[
\exists z\big(z^2=x+y~\wedge~ (u-z)H=(v-z)H\big).
\]
Then we have shown that for $i=1,\ldots,m$ and $S\seq \{1,\ldots,m\}$, $F\models\varphi(a_i;b_S,\sigma(c),c)$ if and only if $a_i+b_S$ is a square in $k$. Since $m\geq 1$ was arbitrary, $\varphi(x;y,u,v)$ has the independence property in $F$, contradicting stability of $F$. 
\end{proof}

\subsection{Minimal bounded fields}
 There is a famous conjecture for minimal fields in a similar spirit to the Stable Fields Conjecture. Recall that a field $F$ is \emph{minimal} if every subset of $F$ definable with parameters in $F$ is finite or cofinite.
\begin{conj}[Podewski~\cite{Pod}]
Every infinite minimal field is algebraically closed.
\end{conj}
Wagner proved the conjecture in positive characteristic~\cite{Wagner-podewski}. Krupi\'nski, Tanovi\'c, and Wagner showed that any counterexample must admit a definable partial order on the field with an infinite chain~\cite[Theorem~1 and Lemma~1.1]{KTW-podewski}. The conjecture
remains open in characteristic zero.

In this section, we apply the ideas in Section~\ref{sec:SFCpos} to study the Podewski conjecture. Recall that a field is \emph{large} if every smooth geometrically integral curve over it with a rational point has infinitely many rational points. Recall that $F$ is \emph{bounded} if it has only finitely many separable extensions of each degree $d$. Koenigsmann showed that every large minimal field is algebraically closed~\cite[\S6.2]{BarySoroker-Fehm}. He also conjectured that every infinite bounded field is large~\cite[p.~496]{Junker-Koenigsmann}. This conjecture remains wide open. If true, it would reduce the bounded case of Podewski's conjecture to the large case. We will give a proof of Podewski's conjecture for bounded fields. In light of Wagner's result in the positive characteristic case, we will restrict to characteristic $0$ for simplicity (however, the proof does not technically require this restriction; see Remark \ref{rem:Podpos}). 

We will need the following facts concerning minimal fields. 

\begin{fact} \label{fact:minimal}
Let $F$ be an infinite minimal field of characteristic $0$.
\begin{enumerate}
    \item $F^\times$ is divisible.   
    \item For all $d\geq 1$, $F$ contains the set $\mu_d$ of $d^{\text{th}}$ roots of unity.
\end{enumerate}
\end{fact} 
\begin{proof}[Explanation]
See \cite[Lemma~3.2]{KTW-podewski} for part $(1)$. By a result of Milliet \cite{Milliet2013}, part $(2)$  holds for the more general class of infinite ``weakly small" fields. But here is a quick argument in our setting via Kummer theory (similar to the proof of Lemma \ref{lem:divisible}).

Suppose $d$ is least such that $\mu_d\not\subseteq F$. Then
$F(\mu_d)/F$ is a nontrivial abelian extension of degree at most
$\varphi(d)<d$. It therefore has a cyclic subextension $L/F$ of prime
degree $\ell<d$. By minimality of $d$, we have $\mu_\ell\subseteq F$,
so Kummer theory gives $L=F(\sqrt[\ell]{a})$ for some $a\in F^\times$. Since $F^\times$ is divisible, this yields $L=F$, contradicting $[L:F]=\ell$.
\end{proof}

We next prove a technical lemma for minimal fields.

\begin{lemma}\label{lem:d-root}
Let $F$ be a minimal field of characteristic $0$ and let $m\in F[X]$ be a monic irreducible polynomial of degree $d>1$. Then there is no definable function $g\colon D\to F$ on a cofinite set $D\subseteq F$ with $g(t)^d=m(t)$ for all $t\in D$.
\end{lemma}
\begin{proof}
Let $C=\{(x,y)\in F^2:y^d=m(x)\}$ and define $\pi\colon C\to F$ so that $\pi(x,y)=x-y$. We claim that each fiber of $\pi$ has size at most $d-1$. Indeed, given $u\in F$, the fiber of $\pi$ above $u$ consists of points $(x,x-u)$ where $x$ is a root of the polynomial
\[
 P_u(X)=m(X)-(X-u)^d.
\]
The leading terms of $m(X)$ and $(X-u)^d$ cancel since these polynomials are both monic of degree $d$. So $P_u$ has degree at most $d-1$. Also, $P_u$ is nonzero since $P_u(u)=m(u)\ne0$ (recall $m$ has no root in $F$). So $P_u$ has at most $d-1$ roots, as desired.

Now, toward a contradiction, suppose there is a definable function $g\colon D\to F$ as in the lemma. By Fact \ref{fact:minimal}(2), we have $\mu_d\seq F$. By assumption, if $\xi\in\mu_d$ and $t\in D$ then $(t,\xi g(t))\in C$.  For $\xi\in\mu_d$, let $\Gamma_\xi$ be the graph of $\xi g$. Note  that $g(t)\neq 0$ for all $t\in D$ since $m$ has no roots in $F$. Therefore, the sets $\Gamma_\xi$ for $\xi\in \mu_d$ are $d$ pairwise disjoint definable subsets of $C$. Moreover,  since $D$ is infinite, each set $\Gamma_\xi$ is infinite and hence $\pi(\Gamma_\xi)$ is infinite since $\pi$ has finite fibers. So each image $\pi(\Gamma_\xi)$ is cofinite. But then the intersection $\bigcap_{\xi\in\mu_d}\pi(\Gamma_\xi)$ is nonempty, and any point in this intersection yields a fiber of $\pi$ containing $d$ distinct points, a contradiction. 
\end{proof}
The main theorem in this section is the following. 

\begin{theorem}\label{thm:minimal-general}
Let $F$ be a minimal field of characteristic $0$, and let $K=F(\alpha)$ be a proper finite extension of degree $d$. Then the map
\[
 F\rightarrow K^\times/({K}^\times)^{d},\qquad t\mapsto (t-\alpha)(K^\times)^d
 \]
has finite fibers. In particular, for some prime $\ell$ dividing $d$, the field $K$ has $|F|$ distinct cyclic extensions of degree $\ell$, hence $F$ is not bounded.
\end{theorem}

\begin{proof}
As in Proposition~\ref{prop:pick} and Lemma~\ref{lem:divisible}, we consider the relation $(x-z)H=(y-z)H$ for a power subgroup $H$. Here we work in $K$, with $z=\alpha$, $H=({K}^\times)^{d}$, and $x,y$ ranging over $F$.

Let $m$ be the monic minimal polynomial of $\alpha$ over $F$. For
$a\in F$, define
\[
 D_a=\left\{t\in F :
 \frac{t-\alpha}{a-\alpha}\in({K^\times)}^{d}\right\}, 
\]
which is definable in $F$ by interpreting $K$ in $F$. Note also that the sets $D_a$ for $a\in F$ are the fibers of the map from the statement of the theorem. We fix $a\in F$ and show that $D_a$ is finite. 

Toward a contradiction, suppose $D_a$ is infinite, hence cofinite. We use $\N$ to denote the norm function for the field extension $K/F$. Note that if $\xi$ is a $d$-th root of unity and $z\in K^\times$ then $\N(\xi z)=\xi^d\N(z)=\N(z)$ since $[K:F]=d$. In particular, we have a well-defined function $h\colon D_a\to F$ such that for $t\in D_a$,  $h(t)=\N(z)$ where $z\in K^\times$ satisfies  $z^d=(t-\alpha)/(a-\alpha)$. Also, $h$ is definable via the determinant formula for the norm. Note further that
$\N(t-\alpha)=m(t)$ for any $t\in F$, and thus if $t\in D_a$ then
\[
 h(t)^d=\N\!\left(\frac{t-\alpha}{a-\alpha}\right)
       =\frac{m(t)}{m(a)}.
\]
By the divisibility of $F^\times$ (Fact \ref{fact:minimal}(1)), there is $c\in F^\times$ with $c^d=m(a)$. Now the definable function $g=ch\colon D_a\to F$ satisfies $g(t)^d=m(t)$ for all $t\in D_a$,  contradicting Lemma~\ref{lem:d-root}. 

This establishes the main statement of the theorem. Since $F$ is infinite, we conclude that $[K^\times:(K^\times)^d]=|F|$, and thus there is some prime $\ell$ dividing $d$ such that $[K^\times:(K^\times)^\ell]=|F|$. In particular, $K^\times/(K^\times)^\ell$ is an $\FFF_\ell$-vector space  of dimension $|F|$, and hence contains $|F|$ distinct $1$-dimensional subspaces. By Kummer theory (see, e.g., \cite[VIII Theorems 8.1,8.2]{LangAlgebra}), it follows that $K$ has $|F|$ distinct cyclic extensions of degree $\ell$. Each of these yields a distinct  extension of $F$ of degree $d\ell$, hence $F$ is not bounded.
\end{proof}

Together with Wagner's result in positive characteristic, Theorem~\ref{thm:minimal-general} gives Podewski's conjecture for bounded fields.

\begin{cor}\label{cor:Pod}
    Every infinite bounded minimal field is algebraically closed.
\end{cor}

\begin{remark}\label{rem:Podpos}
In order to extend the above argument to positive characteristic, one just needs an appropriate modification of Fact \ref{fact:minimal}. More specifically, suppose $F$ is an infinite minimal field of characteristic $p>0$. Then $F$ is perfect by  \cite[Lemma~3.2]{KTW-podewski}, so all finite extensions are separable. Also, it follows from a result of Whaples~\cite[Theorem~5]{Kulhmann-additive} that no finite extension of $F$ has  degree divisible by $p$. Thus the same Kummer argument shows that if $F$ has a finite extension of degree $d$ then $\mu_d\seq F$. This allows one to remove the characteristic $0$ assumption from both Lemma \ref{lem:d-root} and Theorem \ref{thm:minimal-general} (assuming $F$ is infinite).

That said, Wagner's proof in the positive characteristic case is much simpler.
\end{remark}

\section{Final remarks}

Obviously, the main question still open is whether the Stable Fields Conjecture holds in characteristic $0$. 
Note that Theorem \ref{thm:SFCpos} implies that a counterexample must satisfy an axiom scheme (namely, ``characteristic $0$") that fails in every field of positive characteristic. Theorem \ref{thm:Lef} shows that, in fact, a counterexample must satisfy a single sentence (possibly depending on the field) that fails in all positive characteristic fields.

Turning now to $\CCC(t)$, a natural question is whether $\CCC(t)$ is NIP. A negative answer would follow from the NIP Fields Conjecture \cite[Conjecture 1.1]{AnscombeJahnkeNIP}. It is known, however, that  $\CCC(t)$ is not strongly NIP (or even strongly NTP$_2$) due to a general result of Chernikov, Kaplan, and Simon \cite[Corollary 4.6]{CKSNTP2}. We also record the probably folklore observation that $\CCC(t)$ is not distal. Indeed, a transcendence basis for $\CCC$ over $\QQQ$ is an indiscernible set since any permutation extends to an automorphism of $\CCC(t)$, hence distality fails by \cite[Example 2.4]{SimonDistal}. A positive answer to a  general conjecture of Aschenbrenner, Chernikov, Gehret, and Ziegler \cite[Conjecture 6.17]{ACGZDistality} would imply that, in fact, $\CCC(t)$ has no distal expansion. 

\bibliographystyle{plain}
\bibliography{research}

@book{LangAlgebra,
  author    = {Lang, Serge},
  title     = {Algebra},
  edition   = {Revised third},
  series    = {Graduate Texts in Mathematics},
  volume    = {211},
  publisher = {Springer-Verlag},
  address   = {New York},
  year      = {2002},
  doi       = {10.1007/978-1-4613-0041-0}
}

@article{AschenbrennerSchoutens2007,
  author  = {Aschenbrenner, Matthias and Schoutens, Hans},
  title   = {{Lefschetz} extensions, tight closure, and big {Cohen-Macaulay} algebras},
  journal = {Israel J. Math.},
  volume  = {161},
  year    = {2007},
  pages   = {221--310},
  doi     = {10.1007/s11856-007-0080-0}
}

@book{LamFirstCourse,
  author    = {Lam, T. Y.},
  title     = {A First Course in Noncommutative Rings},
  edition   = {Second},
  series    = {Graduate Texts in Mathematics},
  volume    = {131},
  publisher = {Springer-Verlag},
  address   = {New York},
  year      = {2001},
  doi       = {10.1007/978-1-4419-8616-0},
  mrnumber  = {1838439}
}

@article {Junker-Koenigsmann,
    AUTHOR = {Junker, Markus and Koenigsmann, Jochen},
     TITLE = {Schlanke {K}\"orper (slim fields)},
   JOURNAL = {J. Symbolic Logic},
  FJOURNAL = {The Journal of Symbolic Logic},
    VOLUME = {75},
      YEAR = {2010},
    NUMBER = {2},
     PAGES = {481--500},
      ISSN = {0022-4812,1943-5886},
   MRCLASS = {03C60},
  MRNUMBER = {2648152},
MRREVIEWER = {Jonathan\ Kirby},
       DOI = {10.2178/jsl/1268917491},
       URL = {https://doi.org/10.2178/jsl/1268917491},
}

@incollection {BarySoroker-Fehm,
    AUTHOR = {Bary-Soroker, Lior and Fehm, Arno},
     TITLE = {Open problems in the theory of ample fields},
 BOOKTITLE = {Geometric and differential {G}alois theories},
    SERIES = {S\'emin. Congr.},
    VOLUME = {27},
     PAGES = {1--11},
 PUBLISHER = {Soc. Math. France, Paris},
      YEAR = {2013},
      ISBN = {978-2-85629-364-5},
   MRCLASS = {12E30 (12F10 12F12)},
  MRNUMBER = {3203546},
MRREVIEWER = {Ivan\ D.\ Chipchakov},
}

@incollection {Kulhmann-additive,
    AUTHOR = {Kuhlmann, Franz-Viktor},
     TITLE = {Additive polynomials and their role in the model theory of
              valued fields},
 BOOKTITLE = {Logic in {T}ehran},
    SERIES = {Lect. Notes Log.},
    VOLUME = {26},
     PAGES = {160--203},
 PUBLISHER = {Assoc. Symbol. Logic, La Jolla, CA},
      YEAR = {2006},
      ISBN = {978-1-56881-296-0; 1-56881-296-5},
   MRCLASS = {03C60 (12L12)},
  MRNUMBER = {2262319},
MRREVIEWER = {Nicolas\ C.\ Guzy},
}

@article {Wagner-podewski,
    AUTHOR = {Wagner, Frank O.},
     TITLE = {Minimal fields},
   JOURNAL = {J. Symbolic Logic},
  FJOURNAL = {The Journal of Symbolic Logic},
    VOLUME = {65},
      YEAR = {2000},
    NUMBER = {4},
     PAGES = {1833--1835},
      ISSN = {0022-4812,1943-5886},
   MRCLASS = {03C60},
  MRNUMBER = {1812183},
MRREVIEWER = {Luc\ B\'elair},
       DOI = {10.2307/2695078},
       URL = {https://doi.org/10.2307/2695078},
}

@article {KTW-podewski,
    AUTHOR = {Krupi\'nski, Krzysztof and Tanovi\'c, Predrag and Wagner,
              Frank O.},
     TITLE = {Around {P}odewski's conjecture},
   JOURNAL = {Fund. Math.},
  FJOURNAL = {Fundamenta Mathematicae},
    VOLUME = {222},
      YEAR = {2013},
    NUMBER = {2},
     PAGES = {175--193},
      ISSN = {0016-2736,1730-6329},
   MRCLASS = {03C60 (12L12 20A15)},
  MRNUMBER = {3085238},
MRREVIEWER = {David\ A.\ Pierce},
       DOI = {10.4064/fm222-2-4},
       URL = {https://doi.org/10.4064/fm222-2-4},
}

@article {Pod,
    AUTHOR = {Podewski, Klaus-Peter},
     TITLE = {Minimale {R}inge},
   JOURNAL = {Math.-Phys. Semesterber.},
  FJOURNAL = {Mathematisch-Physikalische Semesterberichte zur Pflege des
              Zusammenhangs von Schule und Universit\"at},
    VOLUME = {22},
      YEAR = {1975},
    NUMBER = {2},
     PAGES = {193--197},
      ISSN = {0340-4897},
   MRCLASS = {12L99 (02H15)},
  MRNUMBER = {392962},
MRREVIEWER = {H.\ H.\ Storrer},
}

@misc{DembnerFranklin2026,
  author        = {Dembner, Spencer and Franklin, Cole},
  title         = {{$\mathbb{C}(t)$} has the strict order property},
  year          = {2026},
  note          = {{arXiv}:2611.30227 [math.LO]},
  eprint        = {2611.30227},
  archivePrefix = {arXiv},
  primaryClass  = {math.LO},
  doi           = {10.48550/arXiv.2611.30227},
  url           = {https://arxiv.org/abs/2611.30227}
}

@article{ACGZDistality,
  author  = {Aschenbrenner, Matthias and Chernikov, Artem and
             Gehret, Allen and Ziegler, Martin},
  title   = {Distality in valued fields and related structures},
  journal = {Trans. Amer. Math. Soc.},
  volume  = {375},
  number  = {7},
  year    = {2022},
  pages   = {4641--4710},
  doi     = {10.1090/tran/8661}
}

@article{AnscombeBourbaki,
  author  = {Anscombe, Sylvy},
  title   = {{Shelah}'s conjecture and {Johnson}'s theorem (after {Will Johnson})},
  journal = {Ast{\'e}risque},
  number  = {438},
  year    = {2022},
  pages   = {247--279},
  note    = {S{\'e}minaire Bourbaki, vol.~2021/2022, exp.~1186},
  doi     = {10.24033/ast.1187},
  url     = {https://www.bourbaki.fr/TEXTES/Exp1186-Anscombe.pdf}
}

@article{AnscombeFehmFragments,
  author  = {Anscombe, Sylvy and Fehm, Arno},
  title   = {Interpretations of syntactic fragments of theories of fields},
  journal = {Israel J. Math.},
  year    = {2026},
  doi     = {10.1007/s11856-026-2938-z}
}

@article{AnscombeJahnkeNIP,
  author        = {Anscombe, Sylvy and Jahnke, Franziska},
  title         = {Characterizing {NIP} henselian fields},
  journal       = {Journal of the London Mathematical Society},
  year          = {2024},
  volume        = {109},
  pages         = {e12868},
  doi           = {10.1112/jlms.12868}
}

@article{Ax1968,
  author  = {Ax, James},
  title   = {The elementary theory of finite fields},
  journal = {Ann. of Math. (2)},
  volume  = {88},
  number  = {2},
  year    = {1968},
  pages   = {239--271},
  doi     = {10.2307/1970573}
}

@article{CH99,
  title={Model theory of difference fields},
  author={Chatzidakis, Zo{\'e} and Hrushovski, Ehud},
  journal={Transactions of the American Mathematical Society},
  volume={351},
  number={8},
  pages={2997--3071},
  year={1999}
}

@incollection{Chatzidakis1999,
  author    = {Chatzidakis, Zo{\'e}},
  title     = {Simplicity and independence for pseudo-algebraically closed fields},
  booktitle = {Models and Computability},
  editor    = {Cooper, S. Barry and Truss, John K.},
  series    = {London Mathematical Society Lecture Note Series},
  volume    = {259},
  publisher = {Cambridge University Press},
  address   = {Cambridge},
  year      = {1999},
  pages     = {41--61},
  doi       = {10.1017/CBO9780511565670.004}
}

@article{ChatzidakisPillay1998,
  author  = {Chatzidakis, Zo{\'e} and Pillay, Anand},
  title   = {Generic structures and simple theories},
  journal = {Ann. Pure Appl. Logic},
  volume  = {95},
  number  = {1--3},
  year    = {1998},
  pages   = {71--92},
  doi     = {10.1016/S0168-0072(98)00021-9}
}

@article{ChatzidakisVanDenDriesMacintyre1992,
  author  = {Chatzidakis, Zo{\'e} and van den Dries, Lou and Macintyre, Angus},
  title   = {Definable sets over finite fields},
  journal = {J. Reine Angew. Math.},
  volume  = {427},
  year    = {1992},
  pages   = {107--135},
  doi     = {10.1515/crll.1992.427.107}
}

@article{CherlinShelah1980,
  author  = {Cherlin, Gregory L. and Shelah, Saharon},
  title   = {Superstable fields and groups},
  journal = {Ann. Math. Logic},
  volume  = {18},
  number  = {3},
  year    = {1980},
  pages   = {227--270},
  doi     = {10.1016/0003-4843(80)90006-6}
}

@article{Chernikov2014,
  author  = {Chernikov, Artem},
  title   = {Theories without the tree property of the second kind},
  journal = {Ann. Pure Appl. Logic},
  volume  = {165},
  number  = {2},
  year    = {2014},
  pages   = {695--723},
  doi     = {10.1016/j.apal.2013.10.002}
}

@article{CKSNTP2,
  author  = {Chernikov, Artem and Kaplan, Itay and Simon, Pierre},
  title   = {Groups and fields with {$\mathrm{NTP}_2$}},
  journal = {Proc. Amer. Math. Soc.},
  volume  = {143},
  number  = {1},
  year    = {2015},
  pages   = {395--406},
  doi     = {10.1090/S0002-9939-2014-12229-5}
}

@incollection{deligne-rapoport,
  author    = {Deligne, Pierre and Rapoport, Michael},
  title     = {Les sch{\'e}mas de modules de courbes elliptiques},
  booktitle = {Modular Functions of One Variable {II}},
  series    = {Lecture Notes in Mathematics},
  volume    = {349},
  pages     = {143--316},
  publisher = {Springer-Verlag},
  year      = {1973},
  doi       = {10.1007/978-3-540-37855-6_4},
  note      = {\url{https://doi.org/10.1007/978-3-540-37855-6_4}}
}

@article{Douady1960,
  author  = {Douady, Adrien},
  title   = {Cohomologie des groupes compacts totalement discontinus},
  journal = {S\'eminaire Bourbaki},
  volume  = {5},
  year    = {1960},
  pages   = {287--298},
  note    = {Expos\'e no.~189, d\'ecembre 1959},
  url     = {https://numdam.org/item/SB_1958-1960__5__287_0/}
}

@book{DummitFoote2004,
  author    = {Dummit, David S. and Foote, Richard M.},
  title     = {Abstract Algebra},
  edition   = {Third},
  publisher = {John Wiley \& Sons},
  year      = {2004}
}

@phdthesis{Dupont2015,
  author = {Dupont, Katharina},
  title  = {Definable valuations on {NIP} fields},
  school = {Universit{\"a}t Konstanz},
  year   = {2015},
  url    = {https://d-nb.info/1093432799/34}
}

@incollection{Dupont2017,
  author    = {Dupont, Katharina},
  title     = {Definable valuations induced by definable subgroups},
  booktitle = {Groups, Modules, and Model Theory---Surveys and Recent
               Developments: In Memory of {R{\"u}diger G{\"o}bel}},
  editor    = {Droste, Manfred and Fuchs, L{\'a}szl{\'o} and
               Goldsmith, Brendan and Str{\"u}ngmann, Lutz},
  publisher = {Springer},
  address   = {Cham},
  year      = {2017},
  pages     = {83--108},
  doi       = {10.1007/978-3-319-51718-6_5}
}

@article{DupontHassonKuhlmann2019,
  author        = {Dupont, Katharina and Hasson, Assaf and Kuhlmann, Salma},
  title         = {Definable valuations induced by multiplicative subgroups and {NIP} fields},
  journal       = {Arch. Math. Logic},
  volume        = {58},
  number        = {7-8},
  year          = {2019},
  pages         = {819--839},
  doi           = {10.1007/s00153-019-00661-2},
  eprint        = {1704.02910},
  archivePrefix = {arXiv},
  primaryClass  = {math.LO}
}

@incollection{Duret1980,
  author    = {Duret, Jean-Louis},
  title     = {Les corps faiblement alg{\'e}briquement clos non s{\'e}parablement clos ont la propri{\'e}t{\'e} d'ind{\'e}pendance},
  booktitle = {Model Theory of Algebra and Arithmetic},
  editor    = {Pacholski, L. and Wierzejewski, J. and Wilkie, A. J.},
  series    = {Lecture Notes in Mathematics},
  volume    = {834},
  publisher = {Springer-Verlag},
  address   = {Berlin},
  year      = {1980},
  pages     = {136--162},
  doi       = {10.1007/BFb0090164}
}

@misc{FornasieroKaplanMatthews2026,
  author        = {Fornasiero, Antongiulio and Kaplan, Elliot and Matthews, Angus},
  title         = {Geometric fields, ranks, and generic derivations},
  year          = {2026},
  howpublished  = {arXiv:2605.22725 [math.LO]},
  eprint        = {2605.22725},
  archivePrefix = {arXiv},
  primaryClass  = {math.LO},
  url           = {https://arxiv.org/abs/2605.22725}
}

@article{HaleviHassonStronglyDependent,
  author        = {Halevi, Yatir and Hasson, Assaf},
  title         = {Strongly dependent ordered abelian groups and {Henselian} fields},
  journal       = {Israel Journal of Mathematics},
  year          = {2019},
  volume        = {232},
  pages         = {719--758},
  doi           = {10.1007/s11856-019-1885-3}
}

@article{HaleviPalacin2019,
  author  = {Halevi, Yatir and Palac{\'i}n, Daniel},
  title   = {The dp-rank of abelian groups},
  journal = {J. Symbolic Logic},
  volume  = {84},
  number  = {3},
  pages   = {957--986},
  year    = {2019},
  doi     = {10.1017/jsl.2018.89}
}

@article{Hempel2016,
  author  = {Hempel, Nadja},
  title   = {On {$n$}-dependent groups and fields},
  journal = {Math. Logic Quart.},
  volume  = {62},
  number  = {3},
  year    = {2016},
  pages   = {215--224},
  doi     = {10.1002/malq.201400080}
}

@article{HHJStronglyDependentFields,
  author        = {Halevi, Yatir and Hasson, Assaf and Jahnke, Franziska},
  title         = {A conjectural classification of strongly dependent fields},
  journal       = {Bulletin of Symbolic Logic},
  year          = {2019},
  volume        = {25},
  pages         = {182--195},
  doi           = {10.1017/bsl.2019.13},
  url           = {https://arxiv.org/abs/1805.03814}
}

@article{HHJVTopologies,
  author        = {Halevi, Yatir and Hasson, Assaf and Jahnke, Franziska},
  title         = {Definable {$V$}-topologies, henselianity and {NIP}},
  journal       = {Journal of Mathematical Logic},
  year          = {2020},
  volume        = {20},
  pages         = {2050008},
  doi           = {10.1142/S0219061320500087},
  url           = {https://arxiv.org/abs/1901.05920}
}

@incollection{Hrushovski2002,
  author    = {Hrushovski, Ehud},
  title     = {Pseudo-finite fields and related structures},
  booktitle = {Model Theory and Applications},
  editor    = {B{\'e}lair, Luc and Chatzidakis, Zo{\'e} and D'Aquino, Paola and Marker, David and Otero, Margarita and Point, Fran{\c{c}}oise and Wilkie, A. J.},
  series    = {Quaderni di Matematica},
  volume    = {11},
  publisher = {Aracne},
  address   = {Rome},
  year      = {2002},
  pages     = {151--212},
  mrnumber  = {2159717}
}

@article{IzquierdoLucchiniArteche2025,
  author  = {Izquierdo, Diego and Lucchini Arteche, Giancarlo},
  title   = {Transfer principles for {Galois} cohomology and {Serre}'s conjecture {II}},
  journal = {Adv. Math.},
  volume  = {480},
  year    = {2025},
  pages   = {Paper No.~110532},
  doi     = {10.1016/j.aim.2025.110532}
}

@article{JahnkeHenselianExpansions,
  author        = {Jahnke, Franziska},
  title         = {Henselian expansions of {NIP} fields},
  journal       = {Journal of Mathematical Logic},
  year          = {2024},
  volume        = {24},
  pages         = {2350006},
  doi           = {10.1142/S021906132350006X},
  url           = {https://arxiv.org/abs/1607.02953}
}

@article{JahnkeSimonNIP,
  author        = {Jahnke, Franziska and Simon, Pierre},
  title         = {{NIP} henselian valued fields},
  journal       = {Archive for Mathematical Logic},
  year          = {2020},
  volume        = {59},
  pages         = {167--178},
  doi           = {10.1007/s00153-019-00685-8}
}

@article{JohnsonCanonicalDpMinimal,
  author        = {Johnson, Will},
  title         = {The canonical topology on dp-minimal fields},
  journal       = {Journal of Mathematical Logic},
  year          = {2018},
  volume        = {18},
  pages         = {1850007},
  doi           = {10.1142/S0219061318500071}
}

@article{JohnsonDpFiniteIA,
  author        = {Johnson, Will},
  title         = {Dp-finite fields {I(A)}: The infinitesimals},
  journal       = {Annals of Pure and Applied Logic},
  year          = {2021},
  volume        = {172},
  pages         = {102947},
  doi           = {10.1016/j.apal.2021.102947}
}

@article{JohnsonDpFiniteIB,
  author        = {Johnson, Will},
  title         = {Dp-finite fields {I(B)}: Positive characteristic},
  journal       = {Annals of Pure and Applied Logic},
  year          = {2021},
  volume        = {172},
  pages         = {102949},
  doi           = {10.1016/j.apal.2021.102949}
}

@unpublished{JohnsonDpFiniteVI,
  author        = {Johnson, Will},
  title         = {Dp-finite fields {VI}: The dp-finite {Shelah} conjecture},
  year          = {2020},
  note          = {arXiv:2005.13989},
  url           = {https://arxiv.org/abs/2005.13989}
}

@article{JohnsonDpMinimal,
  author        = {Johnson, Will},
  title         = {The classification of dp-minimal and dp-small fields},
  journal       = {Journal of the European Mathematical Society},
  year          = {2023},
  volume        = {25},
  pages         = {467--513},
  doi           = {10.4171/JEMS/1187}
}

@inproceedings{JohnsonICM2026,
  author    = {Johnson, Will},
  title     = {Henselianity and {NIP} fields},
  booktitle = {Proceedings of the International Congress of Mathematicians 2026},
  editor    = {Friedlander, Susan and Tschinkel, Yuri},
  volume    = {3},
  publisher = {Society for Industrial and Applied Mathematics},
  year      = {2026},
  pages     = {97--112},
  doi       = {10.1137/25M1805412}
}

@misc{JohnsonYe2026,
  author        = {Johnson, Will and Ye, Jinhe},
  title         = {Strongly stable fields},
  year          = {2026},
  howpublished  = {arXiv:2610.41430 [math.LO]},
  eprint        = {2610.41430},
  archivePrefix = {arXiv},
  primaryClass  = {math.LO},
  url           = {https://arxiv.org/abs/2610.41430}
}

@article{JSWDpMinimal,
  author        = {Jahnke, Franziska and Simon, Pierre and Walsberg, Erik},
  title         = {Dp-minimal valued fields},
  journal       = {Journal of Symbolic Logic},
  year          = {2017},
  volume        = {82},
  pages         = {151--165},
  doi           = {10.1017/jsl.2016.15}
}

@article{JTWY2024,
  author  = {Johnson, Will and Tran, Chieu-Minh and Walsberg, Erik and Ye, Jinhe},
  title   = {The {\'e}tale-open topology and the stable fields conjecture},
  journal = {J. Eur. Math. Soc. (JEMS)},
  volume  = {26},
  number  = {10},
  pages   = {4033--4070},
  year    = {2024},
  doi     = {10.4171/JEMS/1345}
}

@misc{KaplanMatthewsWalsberg2026,
  author        = {Kaplan, Elliot and Matthews, Angus and Walsberg, Erik},
  title         = {Building trees in large fields},
  year          = {2026},
  howpublished  = {arXiv:2607.16942 [math.LO]},
  eprint        = {2607.16942},
  archivePrefix = {arXiv},
  primaryClass  = {math.LO},
  url           = {https://arxiv.org/abs/2607.16942}
}

@article{KaplanScanlonWagner2011,
  author  = {Kaplan, Itay and Scanlon, Thomas and Wagner, Frank O.},
  title   = {{Artin--Schreier} extensions in {NIP} and simple fields},
  journal = {Israel Journal of Mathematics},
  volume  = {185},
  number  = {1},
  year    = {2011},
  pages   = {141--153},
  doi     = {10.1007/s11856-011-0104-7}
}

@article{KimPillay1997,
  author  = {Kim, Byunghan and Pillay, Anand},
  title   = {Simple theories},
  journal = {Ann. Pure Appl. Logic},
  volume  = {88},
  number  = {2--3},
  year    = {1997},
  pages   = {149--164},
  doi     = {10.1016/S0168-0072(97)00019-5}
}

@article{KimPillay1998,
  author  = {Kim, Byunghan and Pillay, Anand},
  title   = {From stability to simplicity},
  journal = {Bull. Symbolic Logic},
  volume  = {4},
  number  = {1},
  year    = {1998},
  pages   = {17--36},
  doi     = {10.2307/421004}
}

@article{KKLAlmostRealClosed,
  author        = {Krapp, Lothar Sebastian and Kuhlmann, Salma and Leh{\'e}ricy, Gabriel},
  title         = {Strongly {NIP} almost real closed fields},
  journal       = {Mathematical Logic Quarterly},
  year          = {2021},
  volume        = {67},
  pages         = {321--328},
  doi           = {10.1002/malq.202000060}
}

@article{Krupinski2015,
  author  = {Krupi{\'n}ski, Krzysztof},
  title   = {Superrosy fields and valuations},
  journal = {Ann. Pure Appl. Logic},
  volume  = {166},
  number  = {3},
  year    = {2015},
  pages   = {342--357},
  doi     = {10.1016/j.apal.2014.11.006}
}

@article{KrupinskiPillay2011,
  author  = {Krupi{\'n}ski, Krzysztof and Pillay, Anand},
  title   = {On stable fields and weight},
  journal = {J. Inst. Math. Jussieu},
  volume  = {10},
  number  = {2},
  year    = {2011},
  pages   = {349--358},
  doi     = {10.1017/S1474748010000228}
}

@article{Macintyre1971,
  author  = {Macintyre, Angus},
  title   = {On {$\omega_1$}-categorical theories of fields},
  journal = {Fund. Math.},
  volume  = {71},
  number  = {1},
  year    = {1971},
  pages   = {1--25},
  doi     = {10.4064/fm-71-1-1-25}
}

@article{MartinPizarro2006,
  author  = {Martin-Pizarro, Amador},
  title   = {Elliptic and hyperelliptic curves over supersimple fields in characteristic 2},
  journal = {J. Pure Appl. Algebra},
  volume  = {204},
  number  = {2},
  year    = {2006},
  pages   = {368--379},
  doi     = {10.1016/j.jpaa.2005.05.002}
}

@article{MartinPizarroPillay2004,
  author  = {Martin-Pizarro, Amador and Pillay, Anand},
  title   = {Elliptic and hyperelliptic curves over supersimple fields},
  journal = {J. London Math. Soc. (2)},
  volume  = {69},
  number  = {1},
  year    = {2004},
  pages   = {1--13},
  doi     = {10.1112/S002461070300499X}
}

@article{MartinPizarroWagner2009,
  author  = {Martin-Pizarro, Amador and Wagner, Frank O.},
  title   = {Supersimplicity and quadratic extensions},
  journal = {Arch. Math. Logic},
  volume  = {48},
  number  = {1},
  year    = {2009},
  pages   = {55--61},
  doi     = {10.1007/s00153-008-0115-1}
}

@book{Matsumura1986,
  author    = {Matsumura, Hideyuki},
  title     = {Commutative Ring Theory},
  series    = {Cambridge Studies in Advanced Mathematics},
  volume    = {8},
  note      = {Translated by M. Reid},
  publisher = {Cambridge University Press},
  year      = {1986}
}

@article{Milliet2013,
  author        = {Milliet, C{\'e}dric},
  title         = {Fields with few types},
  journal       = {J. Symbolic Logic},
  volume        = {78},
  number        = {1},
  year          = {2013},
  pages         = {72--84},
  doi           = {10.2178/jsl.7801050},
  eprint        = {1104.0288},
  archivePrefix = {arXiv},
  primaryClass  = {math.LO}
}

@misc{Moerman2024,
  author       = {Moerman, Boaz},
  title        = {Generalized {Campana} points and adelic approximation on toric varieties},
  year         = {2024},
  howpublished = {arXiv:2407.03048},
  url          = {https://arxiv.org/abs/2407.03048}
}

@article{Morley1965,
  author  = {Morley, Michael},
  title   = {Categoricity in power},
  journal = {Trans. Amer. Math. Soc.},
  volume  = {114},
  number  = {2},
  year    = {1965},
  pages   = {514--538},
  doi     = {10.1090/S0002-9947-1965-0175782-0}
}

@article{Nguyen2024,
  author  = {Nguyen, Dong Quan Ngoc},
  title   = {Higher reciprocity law and an analogue of the {Grunwald--Wang} theorem for the ring of polynomials over an ultra-finite field},
  journal = {Annals of Pure and Applied Logic},
  volume  = {175},
  number  = {6},
  year    = {2024},
  pages   = {103438},
  doi     = {10.1016/j.apal.2024.103438}
}

@incollection{Pillay1998,
  author    = {Pillay, Anand},
  title     = {The model-theoretic content of {Lang}'s conjecture},
  editor    = {Bouscaren, Elisabeth},
  booktitle = {Model Theory and Algebraic Geometry: An Introduction to
               {E. Hrushovski}'s Proof of the Geometric {Mordell--Lang} Conjecture},
  series    = {Lecture Notes in Mathematics},
  volume    = {1696},
  publisher = {Springer},
  address   = {Berlin, Heidelberg},
  year      = {1998},
  pages     = {101--106},
  doi       = {10.1007/978-3-540-68521-0_6}
}

@article{PillayPoizat1995,
  author  = {Pillay, Anand and Poizat, Bruno},
  title   = {Corps et chirurgie},
  journal = {J. Symbolic Logic},
  volume  = {60},
  number  = {2},
  year    = {1995},
  pages   = {528--533},
  doi     = {10.2307/2275848}
}

@article{PillayScanlonWagner1998,
  author  = {Pillay, Anand and Scanlon, Thomas and Wagner, Frank O.},
  title   = {Supersimple fields and division rings},
  journal = {Math. Res. Lett.},
  volume  = {5},
  number  = {4},
  year    = {1998},
  pages   = {473--483},
  doi     = {10.4310/MRL.1998.v5.n4.a5}
}

@article{PillayWalsberg2023,
  author  = {Pillay, Anand and Walsberg, Erik},
  title   = {{Galois} groups of large simple fields},
  journal = {Model Theory},
  volume  = {2},
  number  = {2},
  year    = {2023},
  pages   = {357--380},
  note    = {With an appendix by Philip Dittmann},
  doi     = {10.2140/mt.2023.2.357}
}

@book{Poi01,
author = {Poizat, Bruno},
address = {Providence, Rhode Island},
booktitle = {Stable groups},
isbn = {1-4704-1314-0},
language = {eng;fre},
publisher = {American Mathematical Society},
series = {Mathematical surveys and monographs, volume 87},
title = {Stable groups },
year = {2001},
}

@article{Poizat1983,
  author  = {Poizat, Bruno},
  title   = {Groupes stables, avec types g{\'e}n{\'e}riques r{\'e}guliers},
  journal = {The Journal of Symbolic Logic},
  volume  = {48},
  number  = {2},
  year    = {1983},
  pages   = {339--355},
  doi     = {10.2307/2273551}
}

@book{Ramanathan1954,
  author    = {Ramanathan, K. G.},
  title     = {Lectures on the Algebraic Theory of Fields},
  series    = {Lectures on Mathematics and Physics, Mathematics},
  volume    = {5},
  publisher = {Tata Institute of Fundamental Research},
  address   = {Bombay},
  year      = {1954},
  url       = {https://mathweb.tifr.res.in/sites/default/files/publications/ln/tifr05.pdf}
}

@article{scanlon-ct,
  author  = {Scanlon, Thomas},
  title   = {Decidability of some complicated structures definable in {$\mathbb{C}(t)$}},
  journal = {Journal of Symbolic Logic},
  year    = {2026},
  pages   = {1--16},
  doi     = {10.1017/jsl.2026.10279},
  note    = {Published online 30 October 2026. \url{https://doi.org/10.1017/jsl.2026.10279}}
}

@unpublished{Scanlon1999,
  author = {Scanlon, Thomas},
  title  = {Infinite stable fields are {Artin--Schreier} closed},
  year   = {1999},
  note   = {Unpublished note, October~1},
  url    = {https://math.berkeley.edu/~scanlon/papers/ASclos.pdf}
}

@misc{Scanlon2007slides,
  author = {Scanlon, Thomas},
  title  = {Definability in fields. {Lecture} 3: {Finding} structure through transcendence},
  year   = {2007},
  note   = {Lecture slides, Model Theory and Computable Model Theory,
            University of Florida, Gainesville, February 6},
  url    = {https://math.berkeley.edu/~scanlon/papers/dedfg3.pdf}
}

@book{Serre2002,
  author    = {Serre, Jean-Pierre},
  title     = {{Galois} Cohomology},
  series    = {Springer Monographs in Mathematics},
  publisher = {Springer-Verlag},
  address   = {Berlin},
  year      = {2002},
  note      = {Translated from the French by Patrick Ion.
               Corrected second printing},
  doi       = {10.1007/978-3-642-59141-9}
}

@article{Shelah1975,
  author   = {Shelah, Saharon},
  title    = {The lazy model-theoretician's guide to stability},
  journal  = {Logique et Analyse (N.S.)},
  year     = {1975},
  volume   = {18},
  number   = {71--72},
  pages    = {241--308},
  mrnumber = {0539969},
  url      = {https://shelah.logic.at/papers/54/}
}

@article{ShelahStronglyDependent,
  author        = {Shelah, Saharon},
  title         = {Strongly dependent theories},
  journal       = {Israel Journal of Mathematics},
  year          = {2014},
  volume        = {204},
  pages         = {1--83},
  doi           = {10.1007/s11856-014-1111-2}
}

@book{silverman-aec,
  author    = {Silverman, Joseph H.},
  title     = {The Arithmetic of Elliptic Curves},
  edition   = {Second},
  series    = {Graduate Texts in Mathematics},
  volume    = {106},
  publisher = {Springer},
  address   = {New York},
  year      = {2009},
  doi       = {10.1007/978-0-387-09494-6},
  note      = {\url{https://doi.org/10.1007/978-0-387-09494-6}}
}

@article{SimonDistal,
  author  = {Simon, Pierre},
  title   = {Distal and non-distal {NIP} theories},
  journal = {Ann. Pure Appl. Logic},
  volume  = {164},
  number  = {3},
  year    = {2013},
  pages   = {294--318},
  doi     = {10.1016/j.apal.2012.10.015}
}

@phdthesis{Sinclair2018,
  author = {Sinclair, Peter},
  title  = {Relationships between model theory and valuations of fields},
  school = {McMaster University},
  year   = {2018},
  url    = {https://hdl.handle.net/11375/23326}
}

@misc{stacks-project,
  author       = {The {Stacks Project Authors}},
  title        = {{The Stacks Project}},
  howpublished = {\url{https://stacks.math.columbia.edu}},
  year         = {2026},
}

@article{Wagner1998,
  author  = {Wagner, Frank O.},
  title   = {Small fields},
  journal = {J. Symbolic Logic},
  volume  = {63},
  number  = {3},
  year    = {1998},
  pages   = {995--1002},
  doi     = {10.2307/2586722}
}

@article{Wood1979,
  author  = {Wood, Carol},
  title   = {Notes on the stability of separably closed fields},
  journal = {J. Symbolic Logic},
  volume  = {44},
  number  = {3},
  year    = {1979},
  pages   = {412--416},
  doi     = {10.2307/2273133}
}

@article{Yoneda2006,
  author  = {Yoneda, Ikuo},
  title   = {Facts on {PAC} fields and stable fields},
  journal = {RIMS K{\^o}ky{\^u}roku},
  volume  = {1515},
  year    = {2006},
  pages   = {52--62},
  url     = {https://www.kurims.kyoto-u.ac.jp/~kyodo/kokyuroku/contents/pdf/1515-6.pdf}
}

\end{document}